\documentstyle[12pt]{article}
\newcommand{\Ibb}[1]{ {\rm I\ifmmode\mkern
            -3.6mu\else\kern -.2em\fi#1}}
\newcommand{\ibb}[1]{\leavevmode\hbox{\kern.3em\vrule
     height 1.2ex depth -.3ex width .2pt\kern-.3em\rm#1}}
\newcommand{\Cx}{{\ibb C}}

\newcommand{\Rl}{{\Ibb R}}
\newcommand{\be}{\begin{eqnarray}}
\newcommand{\ee}{\end{eqnarray}}
\newcommand{\bez}{\begin{eqnarray*}}
\newcommand{\eez}{\end{eqnarray*}}
\renewcommand{\O}{\Omega}
\newcommand{\si}{\sigma}
\newcommand{\we}{\wedge}
\newcommand{\so}{\star_\O}
\newcommand{\pa}{\partial}
\begin{document}
\renewcommand{\theequation} {\arabic{section}.\arabic{equation}}
\begin{tabbing}
\hspace*{12cm} \= Rethymnon 18/96  \\
               \> May 1996
\end{tabbing}
\vskip1.cm
{\LARGE \bf On the Uniqueness of the Moyal Structure 
of Phase-Space Functions} 
\vskip1cm

\noindent
{\bf C. Tzanakis}\\ 
University of Crete,
74100 Rethymnon,
Crete, Greece
\vskip.2cm

\noindent
{\bf A. Dimakis}\\
Institut F\"ur Theoretische Physik,
Bunsenstr. 9, 37073 G\"ottingen, 
Germany

\vskip1.5cm

\noindent
{\small {\bf Abstract:}

\noindent
It is shown that the only associative algebras with a trivial center
defined on functions of $\Rl^N$ by an integral kernel are generalized
Moyal algebras, corresponding to some particular operator ordering. 
Similarly, the only such Lie algebras are generalized Moyal or Poisson
Lie algebras. In both cases these structures are isomorphic respectively 
to the Moyal algebra, the Moyal Lie algebra and the Poisson Lie algebra.
These results are independent of $N$, which is necessarily even, and 
generalize previous work on the subject. }

\vskip1.5cm

\section{Introduction}
\setcounter{equation}{0}

It is well-known that quantum operators can be mapped to phase-space 
functions and vice versa, in many different ways, depending on the 
ordering rule chosen for the operator to which the monomial $q^n p^m$
is mapped (see e.g.\ \cite{AgWo}). Originally this was done by the Weyl 
transformation corresponding to a symmetric ordering \cite{We}. The
inverse of this mapping, the Wigner transformation, was originally 
devised in order to formulate quantum expectation values as classical 
averages on the phase space $\Gamma$ of the system under consideration
\cite{Wi}. Since the work of Moyal \cite{Mo}, who showed the relation 
between these two mappings, many other ordering rules have been considered,
corresponding to some generalization of the Wigner transformation. 
A large class of such transformations\footnote{It can be shown that 
any linear transformation of quantum operators to phase-space functions,
which is phase-space translation invariant is the inverse of 
(\ref{I1}). The calculations will not be given here, but see
\cite{KrPo} eq.(15).} is given by the formula (\cite{AgWo}, \cite{Du})
\be 
\O:A(\vec{q},\vec{p}) \to \O(A)=\hat{A}:={1\over(2\pi)^n}
\int d\si \O(\si)\tilde{A}(\si) e^{i\si\cdot \hat{z}}\;,
\label{I1}
\ee
where we use the following notation:  
\bez
\begin{array}{rcl}
\mbox{phase-space coordinates:}& \;\; & (\vec{q},\vec{p})=: z \\
\mbox{corresponding quantum operators:}& & (\hat{q},
\hat{p}=-i\hbar{\pa\over\pa\vec{q}})=:\hat{z}\\
\mbox{quantum mechanical operators:} & & \hat{A}\\
\mbox{phase-space functions:} & & A(\vec{q},\vec{p})\\
\mbox{Fourier transform:} & & \tilde{A}(\vec{\eta},\vec{\xi})=
{1\over(2\pi)^n}\int dz e^{-i\si\cdot z} A(z),\;\; 
\si=(\vec{\eta},\vec{\xi})\;.
\end{array}
\eez
Here and in what follows, we consider an $n$-dimensional configuration
space. Expressions like $\si\cdot z$ denote the scalar product in
$\Rl^{2n}$ and $d\vec{q}$ etc.\ is an abbreviation for an $n$-dimenional
differential.

In (\ref{I1}) the kernel $\O(\si)$ is an entire analytic function of
$\si$. If it has no zeros, then the inverse transformation exists and is
given by (\cite{TzGr} section 3)
\be 
\O^{-1}:\hat{A}\to A(\vec{q},\vec{p})={1\over\pi^n}\int d\vec{q}\,'
d\vec{p}\,'d\vec{t} e^{-\vec{p}\,'\cdot\vec{t}/\mu}\omega(\vec{q}-\vec{q}\,',
\vec{p}-\vec{p}\,')\langle \vec{q}\,'-\vec{t}\,|\hat{A}|\vec{q}\,'+\vec{t}
\,\rangle\;,\label{I2}
\ee
where
\be
\omega(z):= {1\over(2\pi)^n}\int d\si {e^{i\si\cdot z}\over\O(\si)}\;,
\qquad\quad \mu:={i\hbar\over2}\;, \label{I3}
\ee
and we use Dirac's notation for the matrix elements of $\hat{A}$. For
$\O=1$, (\ref{I1}), (\ref{I2}) reduce to the Weyl and Wigner 
transformations respectively, hence for $\O\neq 1$ we may call them 
generalized Weyl and Wigner transformartions. 

Since the work of Moyal, it is well-known that defining for any phase-space 
functions $f,g$, the operation $\so$, by
\be  f\so g := \O^{-1}(\O(f)\O(g))\;, \label{I4} \ee
and
\be 
[f,g]_\O := {1\over2\mu}(f\so g-g\so f)\;, \label{I5}
\ee
we endow phase-space functions equipped with the ordinary vector space 
operations, $F(\Gamma)$, with the structure of an associative, in general 
non-abelian algebra, eq.(\ref{I4}), hence also with a Lie-algebra structure,
eq.(\ref{I5}). For $\O=1$ these operations are called the Moyal, product
and bracket respectively, denoted by $\star$, $[\,,]$.   
In the classical limit $\mu\to 0$, the latter reduces to the Poisson bracket.

It can be shown from (\ref{I1}), (\ref{I2}) that (\cite{TzGr}, section 3):
\be 
(f\so g)(z) & = & {1\over (2\pi)^{2n}} \int d\si d\si' \tilde{f}(\si)
\tilde{g}(\si') B(\si,\si') e^{i(\si+\si')\cdot z}\;,\label{I6}\\
\lbrack f,g\rbrack_\O (z) & = & {1\over (2\pi)^{2n}} \int d\si d\si' \tilde{f}(\si)
\tilde{g}(\si') A(\si,\si') e^{i(\si+\si')\cdot z}\;,\label{I7}
\ee
where
\be
B(\si,\si') & = & {\O(\si)\O(\si')\over\O(\si+\si')}\,e^{\mu\,\si'\we\si}\;,
\label{I8}\\
A(\si,\si') & = & {\O(\si)\O(\si')\over\O(\si+\si')}\, 
{\sinh(\mu\,\si'\we\si)\over\mu}\;,\label{I9}
\ee
with $\si'\we\si:=J_{ij}\si'{}^i\si^j$, where $J_{ij}$ is the canonical
symplectic matrix of $\Rl^{2n}$, that is
$J = \left(\begin{array}{cc} 0 & \delta_{ij} \\ -\delta_{ij} & 0 \end{array}
\right)$. 

These imply that the mapping $U$,
\be
f\to Uf:\qquad (Uf)(z) = \int d\si\,\O(\si)\tilde{f}(\si)\,e^{i\si\cdot z}
\;,\label{I10}
\ee
is an algebra and Lie-algebra isomomorphism of $(F(\Gamma),\so,[\,,]_\O)$
and $(F(\Gamma),\star,[\,,])$.
From this follows that a necessary and sufficient condition for 
$[\,,]_\O$ to reduce to the Poisson bracket, given that
$(\vec{q},\vec{p})\stackrel{\O}{\rightarrow}(\hat{q},\hat{p})$,
is that $\lim_{\mu\to 0}\O(\si)=1$ (\cite{TzGr} propositions 3.1, 3.2).

In general, if $\O(\si)\to\O_0(\si)$ as $\mu\to 0$ then 
\be
\lim_{\mu\to 0} A(\si,\si')=
{\O_0(\si)\O_0(\si')\over\O_0(\si+\si')}\si'\we\si\;,\label{I9'}
\ee 
and (\ref{I10}) with $\O_0$, replacing $\O$, is a Lie algebra isomorphic
with the Poisson Lie algebra. We may call such algebras {\em generalized
Poisson Lie algebras}, and consider them as singular limits of 
generalized Moyal algebras in the sense that they cannot be defined
through an associative product. We agree to include  them in (\ref{I9})
for $\mu=0$, and for the sake of brevity we will use only the term
Moyal Lie algebra.

Therefore, although (\ref{I6}), (\ref{I7}) define binary operations 
of functions in classical phase-space that are very useful in various 
problems of quantum physics, the underlying abstract algebraic structure
is independent of $\O$. For a detailed treatment of algebraic and
topological questions, as well as applications to physical problems
see \cite{BFFLS}, \cite{BFFLSa}. Thus the question naturally arises, 
{\em whether more general binary operations are defined via (\ref{I6}), 
(\ref{I7}), which are respectively an associative product and a Lie product 
and for which the corresponding kernels are not of the form (\ref{I8}),
(\ref{I9})} and if so, to characterize the corresponding algebras.

It is the aim of this paper to study this problem, which is closely related 
to that considered in \cite{Fl} as we will show in the next section.

\section{The problem}    
\setcounter{equation}{0}

The linear, binary operations on $F(\Gamma)$ defined by (\ref{I6}),
(\ref{I7}), where $B,\;A$ are {\em not} a priori given by (\ref{I8}),
(\ref{I9}), will be denoted by $\star_B$, $[\,,]_A$ respectively.

Using (\ref{I6}), (\ref{I7}), a straightforward calculation gives 
\bez
f\!\star_B\!(g\!\star_B\!h)(z) & = & {1\over(2\pi)^{3n}}\int d\si d\si' d\si''
e^{i(\si+\si'+\si'')\cdot z}B(\si'',\si')B(\si''+\si',\si)
\tilde{f}(\si)\tilde{g}(\si')\tilde{h}(\si''),\\
\lbrack f,\lbrack g,h \rbrack_A \rbrack_A(z) & = & {1\over(2\pi)^{3n}}
\int d\si d\si' d\si''
e^{i(\si+\si'+\si'')\cdot z}\,A(\si,\si'+\si'')A(\si',\si'')\,
\tilde{f}(\si)\tilde{g}(\si')\tilde{h}(\si''),
\eez
(cf.\ (5.14), (5.15a) of \cite{TzGr}). Elementary manipulations of these 
expressions imply that associativity and Jacobi's identity for
(\ref{I6}), (\ref{I7}) are respectively equivalent to
\be 
B(\si,\si')B(\si+\si',\si'') =  B(\si,\si'+\si'')B(\si',\si'')\;,\label{P1}
\ee
\be 
A(\si,\si'+\si'')A(\si',\si'')+A(\si',\si''+\si)A(\si'',\si)
+A(\si'',\si+\si')A(\si,\si')=0\;,\label{P2}
\ee
\be   A(\si,\si')=-A(\si',\si)\;.  \label{P3}  \ee
To connect the above results with the previous works (cf.\ especially
\cite{AgWo}, \cite{Fl}), we will asume that {\em $B,\;A$ are entire
analytic functions of their arguments and $B$ has no zeros.} In the case of
a 1-dimensional configuration space, developping $A$ in (\ref{I7}) in
a power series around zero, we get
\be 
[f,g]_A(z) = \sum_{r=0}^{+\infty} \sum_{s=0}^{+\infty}
\sum_{j=0}^r \sum_{k=0}^s b_{rj,sk}
(\pa_q^j\pa_p^{r-j}f)(\pa_q^k\pa_p^{s-k}g)\;,\label{P4}
\ee
\be
b_{rj,sk}=\left({r\atop j}\right)\left({s\atop k}\right)(-i)^{r+s}
\pa_\eta^j\pa_\xi^{r-j}\pa_{\eta'}^k\pa_{\xi'}^{s-k}A\,|_{(0,0,0,0)}\;,
\label{P5}
\ee
and $\pa_q f=\pa f/\pa q$ etc.

Equation (\ref{P4}) is identical with the Lie product of a 2-index
infinite Lie algebra considered in \cite{Fl}, provided that the summation
in $r$ and $s$ starts from $r=s=1$, which is easily seen to be equivalent
to the requirement
\be [f,1]_A = 0\;,\qquad \mbox{for all $f\in F(\Gamma)$}\;.\label{P6} \ee
Conversely, for any such algebra, (\ref{P5}) defines the derivatives at
zero of an entire function $A$, hence (\ref{P4}) reduces to (\ref{I7}),
provided of course that
\bez
\sum_{r=0}^{+\infty} \sum_{s=0}^{+\infty}
\sum_{j=0}^r \sum_{k=0}^s 
b_{rj,sk} \eta^j\xi^{r-j}\eta'{}^k\xi'{}^{s-k}
\eez 
converges absolutely.

In view of the above discussion, the result of \cite{Fl} can be formulated
as follows:
\vskip.1cm

\noindent
{\bf Proposition:}{\em If $\Gamma$ is 2-dimensional, a nontrivial 
$[\,,]_A$-Lie algebra is isomorphic either to the Moyal Lie algebra or the 
Poisson Lie algebra.}
\vskip.1cm

\noindent
The method followed consists essentially in substituting (\ref{P4}) in
(\ref{P2}) and solving the resulting recurrent relations for the 
coefficients $b_{rj,sk}$. It should be stressed however, that the
proof has been explicited only for the case $b_{rj,sk}=\delta_{rs}
b_{rjk}$. In addition this method, besides leading to cumbersome 
calculations, becomes unhopefully complicated when applied to the case of 
a multidimensional configuration space. Here we adopt a different approach,
namely, we try to solve the functional equations (\ref{P1})--(\ref{P3})
without using series expansions.

Specifically, in the next section we show that {\em any 
$\star_B$-algebra with a trivial center, is a generalized Moyal algebra},
and in the last section we show that {\em any $[\,,]_A$-Lie algebra 
with a trivial center is a generalized Moyal Lie algebra}. In all cases
(\ref{I10}) shows that, up to isomorphism the Moyal and Poisson  
structures are unique.

\section{The characterization of the $\star_B$-algebras}
\setcounter{equation}{0}

In this section we will show that {\em (\ref{P1}) essentially implies
(\ref{I8})} --- the converse is trivial. This will be done in several
steps:
\vskip.1cm\noindent
(i) Putting $\si=\si'=0$ and $\si'=0$ in (\ref{P1}) we get $B(0,\si'')
=B(0,0)$, $B(\si,0)=B(0,\si'')$ so that 
\be B(0,\si)=B(\si,0)=B(0,0)\;. \label{B1} \ee
(ii) Making a cyclic permutation of $\si,\,\si',\,\si''$ in (\ref{P1})
and multiplying the resulting three expressions we get 
\be 
\lefteqn{B(\si+\si',\si'')B(\si'+\si'',\si)B(\si''+\si,\si')=}\hspace{2cm}
\nonumber\\
& & B(\si,\si'+\si'')B(\si',\si''+\si)B(\si'',\si+\si')\;.\label{B2}
\ee
(iii) Since $B(\si,\si')$ is an entire function without zeros, we may 
write it as (see e.g.\ \cite{Kn}, theorem 1 in section 1, extended by
induction on the number of variables)
\be 
B(\si,\si')= e^{b(\si,\si')}= e^{b_s(\si,\si')}e^{b_a(\si,\si')}\;,
\label{B3}
\ee
where $b$ is an entire function and $b_s$, $b_a$ its symmetric and 
antisymmetric part respectively. Substituting (\ref{B3}) in (\ref{P1})
we get 
\bez
\lefteqn{b_s(\si,\si')-b_s(\si',\si'')+b_s(\si+\si',\si'')-
b_s(\si,\si'+\si'')=}\hspace{2cm}& &\\
& & b_a(\si',\si'')-b_a(\si,\si')+b_a(\si,\si'+\si'')-b_a(\si+\si',\si'')\;.
\eez
Since the l.h.s.\ and the r.h.s.\ are respectively antisymmetric and 
symmetric in $\si,\,\si''$, both are zero,
\be
b_a(\si,\si')-b_a(\si',\si'')+b_a(\si+\si',\si'')-b_a(\si,\si'+\si'')
& = & 0\;, \label{B4}\\
b_s(\si,\si')-b_s(\si',\si'')+b_s(\si+\si',\si'')-b_s(\si,\si'+\si'')
& = & 0\;. \label{B5}
\ee
(iv) Using (\ref{B3}) in (\ref{B2}) we obtain
\be
b_a(\si+\si',\si'')+b_a(\si'+\si'',\si)+b_a(\si''+\si,\si')=0\;. \label{B6}
\ee
Adding (\ref{B4}), (\ref{B6}) we get
\be
b_a(\si''+\si,\si') = b_a(\si'',\si') + b_a(\si,\si')\;. \label{B7}
\ee
From this and the continuity of $b_a$, using standard arguments we conclude
that $b_a$ is a bilinear, antisymmetric function, i.e.\ it is an exterior
2-form on the even dimensional manifold $\Gamma$. Therefore, {\em if it is
nondegenerate} then there is a canonical basis such that (cf.\ (\ref{I8}))
\be b_a(\si,\si')= \mu \si'\we\si\;,\qquad\quad \mu\in\Cx\;.\label{B8}\ee
(v) Acting on (\ref{B5}) with $\pa_\si^2\pa_{\si'}-\pa_{\si'}^2\pa_\si$
we readily obtain\footnote{Here partial derivations denote $2n$-dimensional 
gradients, i.e.\ (\ref{B9}) is a 3-index tensor relation.}
\be
(\pa_\si^2\pa_{\si'}-\pa_{\si'}^2\pa_\si)b_s(\si,\si')=
(\pa_\si^2\pa_{\si'}-\pa_{\si'}^2\pa_\si)b_s(\si,\si'+\si'')\;.\label{B9}
\ee
Therefore the r.h.s.\ is independent of $\si''$, hence of $\si'$ as well.
Thus
\be 
(\pa_\si^2\pa_{\si'}-\pa_{\si'}^2\pa_\si)b_s(\si,\si')=:f(\si)\;.
\label{B10}
\ee
Interchanging $\si,\,\si'$ we get $f(\si)=-f(\si')$ by the symmetry of $b_s$, 
hence $f(\si)=f(0)=0$. Therefore 
\bez
(\pa_\si-\pa_{\si'})\pa_{\si'}\pa_\si b_s(\si,\si')=0\;, \eez
so that 
\be b_s(\si,\si')= \chi(\si+\si')+g(\si)+h(\si')\;.\label{B11} \ee
By the symmetry of $b_s$, $g(\si)=h(\si)+c$. However (\ref{B5}) with
$\si'+\si''=0$ and with the aid of (\ref{B1}), (\ref{B3}) gives
\bez
b_s(\si+\si',-\si')+b_s(\si,\si')= b_s(0,0)+b_s(\si',-\si')\;. \eez
Substituting (\ref{B11}) to this, we get
\bez
\chi(\si)+\chi(\si+\si')+2\chi(0)=h(\si)+h(\si+\si')+2h(0)\;. \eez
Putting $\si'=0$ we finally obtain $h(\si)+h(0)=\chi(\si)+\chi(0)$, hence
by properly absorbing a constant
\be b_s(\si,\si')=-\chi(\si+\si')+\chi(\si)+\chi(\si')\;.\label{B12}\ee
This together with (\ref{B3}), (\ref{B8}) shows that $B(\si,\si')$
is of the form (\ref{I8}), which was to be proved. When {\em $b_a$ is
degenerate, (\ref{B8}) holds on $\Gamma/{\rm Ker}\,b_a$.} 

Summarizing our results, we have that there exist, not uniquely determined,
coordinates, such that $\si=(\rho,\tau)$, with $(\rho,0)\in{\rm Ker}\,b_a$
and 
\be
B(\si,\si') = {\O(\si)\O(\si')\over\O(\si+\si')}\,
e^{b_a(\si,\si')}\;,\qquad \O(\si)=e^{\chi(\si)}\;,\label{B13}
\ee
where $b_a(\si,\si')=\mu\,\tau'\we\tau$ (cf.\ (\ref{I8})). If accordingly,
the dual splitting of the phase-space coordinates is $z=(x,y)$, 
then it is easily seen via
the Fourier transformation that functions of $x$ only --- in other words
functions which satisfy $b_a^{ij}(\pa_i f)=0$ --- belong to the 
center\footnote{Notice that conversely we have already shown that 
nondegeneracy of $b_a$ implies that the algebra is a generalized Moyal
algebra, hence by (\ref{I4}) it has a trivial center.}
of $(F(\Gamma),\star_B)$. Therefore our results can be restated as
\vskip.1cm

\noindent
{\bf Theorem 1:} {\em Any $\star_B$-associative algebra, having a trivial center,
and for which $B(\si,\si')$ is an entire analytic function without zeros,
is a generalized    
Moyal algebra (\ref{I6}), (\ref{I8}), hence by (\ref{I10}) it is 
isomorphic to the Moyal algebra defined by (\ref{I6}), (\ref{I8}) with
$\O=1$.}
\vskip.2cm

We may notice here that the following considerations yield a deeper 
insight to the $\star_B$ algebras: If $\ast$ is {\em any binary operation}
that makes $F(\Gamma)$ an associative algebra over $\Cx$, then
\bez 
(f\ast g)(z)={1\over(2\pi)^{2n}}\int d\si d\si' \tilde{f}(\si)  
\tilde{g}(\si')\, e^{i\si\cdot z}\ast e^{i\si'\cdot z}\;. 
\eez
A $\star_B$-algebra corresponds to the case
\be 
e^{i\si\cdot z}\ast e^{i\si'\cdot z} = B(\si,\si')\,e^{i(\si+\si')\cdot z}\;.
\label{B14}
\ee
A necessary and sufficient condition for (\ref{B14}) is readily seen to be
that $\pa/\pa z^i$  are derivations of the $\ast$-operation as well. By
(\ref{B14}), this operation is somehow related to the additive group 
structure of $\Rl^{2n}$. In fact if we consider the extension $E$ of the
translation group $(\Rl^{2n},+)$ by the multiplicative group $\Cx^\ast=
\Cx-\{0\}$, i.e.\ a short exact sequence 
\bez    1\to \Cx^\ast\to E\to \Rl^{2n}\to 0\;,  \eez
then by writing $(\si,\zeta)$ for an element of $E$, and defining 
\be 
(\si,\zeta)(\si',\zeta') := (\si+\si',B(\si,\si')\,\zeta\zeta')\;,
\label{B15}
\ee
we can show, that the action of $E$ on $\Cx^\ast$ is trivial, i.e.  
$\si\zeta=\zeta$ and that the associativity of (\ref{B15}) is equivalent
to our conditions (\ref{P1}) (see \cite{Ja} section 6.10, eq.(61)).
Thus, determination of $B$, or for that matter, $b$, is equivalent to the 
determination of all equivalence classes of extensions of $(\Rl^{2n},+)$
by $\Cx^\ast$, which is known to be in 1-1 correspondence with the second 
cohomology group $H^2(\Rl^{2n},\Cx^\ast)$ (\cite{Ja} theorem 6.15). 

Let $\delta$  be the coboundary operator of the complex 
$C(\Rl^{2n},\Cx^\ast)$, i.e.\ of the functions from $(\Rl^{2n})^k$ to  
$\Cx^\ast$, $k=0,1,2,\ldots$. Since the action of $\Rl^{2n}$ on $\Cx^\ast$
is trivial, we immediately get from the definition of $\delta$ that for 
$b\in{\rm Ker}\,\delta^{(2)}$
\be 
(\delta b)(\si,\si',\si'') = b(\si',\si'')-b(\si+\si',\si'')+
b(\si,\si'+\si'')-b(\si,\si')=0\;, \label{B16}
\ee
which in view of (\ref{B3}) is identical to (\ref{P1}). On the other hand,
for any $\chi\in C^1(\Rl^{2n},\Cx^\ast)$ we have 
\be
(\delta\chi)(\si,\si')=\chi(\si')-\chi(\si+\si')+\chi(\si)\;. \label{B17}
\ee
Since $\delta^2=0$, elements $b$ of $H^2(\Rl^{2n},\Cx^\ast)$ are  
determined up to $\delta\chi$, a fact equivalent to the isomorphism 
(\ref{I10}) of generalized Moyal algebras. 

In the light of the above remarks, our result can be restated by saying that 
the equivalence classes of $\star_B$-algebras (cf.\ (\ref{I10})) are in
1-1 correspondence with the elements of $H^2(\Rl^{2n},\Cx^\ast)$
depending on the rank of $b_a$ (cf.\ (\ref{B13})).

As a final remark we notice that for $b_a$ nondegenerate (cf.\ (\ref{B8})),
the extension (\ref{B15}) of $\Rl^{2n}$ by $\Cx^\ast$ is the direct
product of $\Rl^\ast$ with the Heisenberg group, a fact following from
(\ref{B15}) in view (\ref{B13}) (see \cite{GuSt} section 15,
particularly eq.(15.2)).

\section{The characterization of the $[\,,]_A$-Lie algebras}
\setcounter{equation}{0}

We next turn to the study of the Lie algebras defined by (\ref{I7}),
i.e.\ to the study of (\ref{P2}), (\ref{P3}), assuming that constants
annihilate the Lie product, i.e\ (\ref{P6})
holds. We readily see that this is equivalent to 
\be    A(0,\si)=0\;. \label{A1} \ee
(i) Differentiating (\ref{P2}) with respect to $\si^i$ and putting $\si=0$
we get
\bez
[\pa^1_i A(0,\si'+\si'')-\pa^1_i A(0,\si')-\pa^1_i A(0,\si'']
A(\si',\si'')=0\;,
\eez
where in this section we write $\pa^a_iA$ for the the $i$-th component of 
the gradient of $A$ in the $a$-th argument ($a=1,2$). This leads
to\footnote{Notice that $A$ cannot be identically zero in an open region
of $\Rl^{2n}\times\Rl^{2n}$, by its analyticity (see e.g.\ \cite{Fu}
theorem 3.13). Hence the zeros of $A$ are accumulation points of nonzero 
points, and then (\ref{A2}) follows by continuity.}
\be
\pa^1_i A(0,\si'+\si'')=\pa^1_i A(0,\si') + \pa^1_i A(0,\si'')\;, 
\label{A2}
\ee
hence to the linearity of $\pa^1_i A(0,\si)$. Accordingly setting
\be \pa^1_i A(0,\si) =\omega_{ij}\si^j\;, \label{A3} \ee
(the summation convention always assumed) and differentiating (\ref{P3}) 
with respect to $\si^i$, $\si'{}^j$ at $\si=\si'=0$ we find from
(\ref{A3}) that
\be   \omega_{ij}+\omega_{ji}=0\;. \label{A4} \ee
(ii) Differentiating (\ref{P2}) with respect to $\si^i$ and $\si^j$ 
at $\si=0$ we obtain after some reductions that
\be
(X_{ij}(\si')+X_{ij}(\si''))A(\si',\si'')=
(a_{ij}(\si')+a_{ij}(\si'')-a_{ij}(\si'+\si''))A(\si',\si'')\;,\label{A5}
\ee
where
\newcounter{costas}
\setcounter{costas}{\value{equation}}
\renewcommand{\theequation} {\arabic{section}.\arabic{costas}\alph{equation}}
\setcounter{equation}{0}
\be 
a_{ij}(\si) & := & \pa^1_i\pa^1_j A(0,\si)\;, \label{A5a} \\
X_{ij}(\si) & := & \si^k\,\omega_{k(i}\delta^\ell_{j)}
{\pa\over\pa\si^\ell}\;.\label{A5b}
\ee
To simplify the notation by supressing indices whenever it is necessary,
we introduce symmetric parameters $\alpha^{ij}$, $\beta^{ij}$ and set
\addtocounter{costas}{1}
\setcounter{equation}{0}
\be
X_\alpha & := & {1\over2}\alpha^{ij}\,X_{ij}(\si)
=\si^k\omega_{ki}\alpha^{ij}{\pa\over\pa\si^j}\;, \label{A6a}\\
Z_\alpha(\si,\si') & := & X_\alpha(\si) + X_\alpha(\si')\;, \label{A6b}\\
a_\alpha & := & {1\over2}\alpha^{ij}\,a_{ij}(\si)\;, \label{A6c}\\
\hat{a}_\alpha(\si,\si') & = & a_\alpha(\si)+a_\alpha(\si')-
a_\alpha(\si+\si')\;.\label{A6d}
\ee
In this notation, (\ref{A5}) takes the equivalent compact form  
\renewcommand{\theequation} {\arabic{section}.\arabic{equation}}
\setcounter{equation}{\value{costas}}
\be 
(Z_\alpha A)(\si',\si'') = (\hat{a}_\alpha A)(\si',\si'')\;. \label{A7}
\ee
(iii) The crucial next step is to notice that by (\ref{A6a}), (\ref{A5b})
implies 
\addtocounter{equation}{1}
\setcounter{costas}{\value{equation}}
\renewcommand{\theequation} {\arabic{section}.\arabic{costas}\alph{equation}}
\setcounter{equation}{0}
\be [ X_\alpha,X_\beta] =X_\gamma\;, \label{A8a} \ee
hence 
\be [Z_\alpha,Z_\beta] = Z_\gamma\;, \label{A8b} \ee
\renewcommand{\theequation} {\arabic{section}.\arabic{equation}}
\setcounter{equation}{\value{costas}}
with
\be 
\gamma^{ij}= \alpha^{ik}\omega_{k\ell}\beta^{\ell j} -
\beta^{ik}\omega_{k\ell}\alpha^{\ell j}\;,\label{A9}
\ee
which is again symmetric and constant. Here $[\,,]$ denotes the 
commutator of two 
operators. Therefore (\ref{A8a}) says that the $X_{ij}$ generate a 
Lie algebra. Moreover, by (\ref{A4})
\be 
(X_{ij}(\si')+X_{ij}(\si''))\omega_{k\ell}\si'{}^k\si''{}^\ell =0\;.
\label{A10}
\ee
{\em Assuming that $\omega_{ij}$ is nondegenerate}, (i.e.\ by (\ref{A4}) 
it is a symplectic form on $\Rl^{2n}$), (\ref{A10}) implies that this Lie
algebra is a subalgebra of the $(2n^2+n)$-dimensional Lie algebra of the
symplectic group of $\omega_{ij}$. Since this is exactly the number of the
independent {\em generators} $X_{ij}$, the latter {\em span the whole 
symplectic Lie algebra} (For general information on the symplectic 
group see e.g.\ \cite{FuHa}, Part III, ch.16,17).

\noindent
(iv) Since $Z_\alpha,\,Z_\beta$ are derivations, applying 
$[Z_\alpha,Z_\beta]$ to $A(\si',\si'')$ and using (\ref{A8b}), (\ref{A7}),
we get
\be 
Z_\alpha\,\hat{a}_\beta - Z_\beta \,\hat{a}_\alpha = \hat{a}_\gamma\;.
\label{A11}
\ee
But, from (\ref{A5b}), (\ref{A6b}) we have
\bez
Z_\alpha(\si',\si'')\,a_\beta(\si') & = & (X_\alpha\,a_\beta)(\si')\;,\\
Z_\alpha(\si',\si'')\,a_\beta(\si'+\si'') & = & (X_\alpha\,a_\beta)
(\si'+\si'')\;, 
\eez
and consequently (\ref{A11}) is rewritten as
\bez
\lefteqn{(X_\alpha\,a_\beta -X_\beta\,a_\alpha-a_\gamma)(\si')+
(X_\alpha\,a_\beta -X_\beta\,a_\alpha-a_\gamma)(\si'')=}\hspace{2cm}&&\\
&& (X_\alpha\,a_\beta -X_\beta\,a_\alpha-a_\gamma)(\si'+\si'')\;.
\eez
Therefore
\be
(X_\alpha\,a_\beta -X_\beta\,a_\alpha-a_\gamma)(\si)=
\si^i\,c_i(\alpha,\beta)\;. \label{A12}
\ee
Differentiating (\ref{A12}) with respect to $\si^i$ at $\si=0$, we get
\bez
\si^i\,c_i(\alpha,\beta) = (X_\alpha k_\beta -X_\beta k_\alpha -k_\gamma)
(\si)\;,
\eez
where
\bez k_\alpha(\si) := \si^i\pa_i a_\alpha(0)\;. \eez
Setting 
\bez \tilde{a}_\alpha(\si) := a_\alpha(\si) -k_\alpha(\si)\;, \eez
equation (\ref{A12}) becomes 
\be 
X_\alpha\,\tilde{a}_\beta - X_\beta\,\tilde{a}_\alpha = \tilde{a}_\gamma\;.
\label{A13}
\ee
As a consequence of (\ref{A13}), (\ref{A8a}) the system of first order 
differential equations (cf.\ (\ref{A6a}))
\be 
X_{ij}\chi(\si) =\tilde{a}_{ij}(\si)\;, \label{A14}
\ee
is locally integrable. 
Going back to (\ref{A5}) and using (\ref{A14}) we may rewrite it locally 
as
\be
[X_{ij}(\si')+X_{ij}(\si'')]\left({\O(\si'+\si'')\over\O(\si')\O(\si'')} 
A(\si',\si'')\right)=0\;, \label{A15}
\ee
where $\O(\si):=e^{\chi(\si)}$. This means that $\O(\si'+\si'')A(\si',\si'')
/\O(\si')\O(\si'')$ is invariant under the symplectic group of $\omega_{ij}$
(cf.\ (\ref{A10})), and therefore it is a function of $\omega(\si',\si'')$,
since this is the only bilinear invariant of the group on  
$\Rl^{2n}\times\Rl^{2n}$ (cf.\ \cite {FuHa}, Appendix F). 
Therefore
\be 
A(\si,\si') = {\O(\si)\O(\si')\over\O(\si+\si')} h(\omega(\si,\si'))\;,
\label{A16}
\ee
for some function $h$. Evidently $h$ satisfies the Jacobi condition 
(\ref{P2}). Then by the lemma in the appendix,
$h(x)=c\sinh\,\mu x$ or $h(x)=c x$ and therefore locally $A$ has 
the form (\ref{I9}) or (\ref{I9'}), in the sense that for any specified
point of $\Rl^{2n}\times\Rl^{2n}$, $\chi$ in (\ref{A14}) and $c,\,\mu$
above are defined in some neighborhood of it.

To get (\ref{A16}) {\em globally}, we notice that the symplectic group of
$\omega_{ij}$ acts transitively on $\Rl^{2n}-\{0\}$; in other words, 
$X^i_{jk}(\si)$ as a $2n\times(2n^2+n)$-matrix has rank $2n$ for $\si\neq0$,
and consequently constants are the only invariants it has. Therefore if
$\chi_a,\;\chi_b$ are local solutions of (\ref{A14}) in some open subsets
$U_a,\,U_b$ of $\Rl^{2n}$, with $U_a\cap U_b\neq\emptyset$, then by 
(\ref{A14}), in $U_a\cap U_b$
\bez   X_{ij}(\chi_a-\chi_b)=0\;, \eez
hence $\chi_a=\chi_b+\kappa_{ab}$ for some constant $\kappa_{ab}$. 
Now since $\Rl^{2n}$ is simply connected 
we can write in a consistent way $\kappa_{ab}=\kappa_b-\kappa_a$ for all 
pairs of open sets $U_a,\,U_b$  with $U_a\cap U_b\neq\emptyset$, of an open 
cover of $\Rl^{2n}$  and thus redefining 
$\chi_a$ to $\chi_a+\kappa_a=\chi_b+\kappa_b$ we 
obtain a global solution $\chi$ of (\ref{A14}), hence $h$ in (\ref{A16}) is
globally defined as well.

Evidently from our results above, follows that if $(F(\Gamma),[\,,])$
has a nontrivial center, then $\omega$ is degenerate (cf.\ last
footnote of the previous section). Conversely, suppose $\omega$ has a 
$d$-dimensional kernel. Then as in the previous section, we may write
$\si=(\rho,\tau)$ with $(\rho,0)\in{\rm Ker}\,\omega$, and similarly
$z=(x,y)$. Take any function $f$ of $x$, then its Fourier transform is 
$\tilde{f}(\rho)\delta(\tau)$ and a direct calculation shows that
for any $(\rho_0,0)\in{\rm Ker}\,\omega$
\bez
F(x):=\tilde{f}(0)(d-i\rho_0\cdot x)-(\pa_k\tilde{f})(0)\rho_0^k=
{\pa\over\pa\rho^k}[\tilde{f}(\rho)(\rho^k-\rho_0^k)
e^{i\rho\cdot x}]|_{\rho=0}\;,
\eez
is in the center of the Lie algebra.   

Therefore we summarize the results of this section in 
\vskip.1cm

\noindent
{\bf Theorem 2:} {\em A $[\,,]_A$-Lie algebra, for which $A$ is an 
entire analytic function and which has a trivial center, 
is a generalized Moyal Lie algebra (\ref{I9}), (\ref{I9'}) and is 
isomorphic  to the Moyal Lie algebra, via (\ref{I10}).}

\section*{Appendix}
\renewcommand{\theequation} {A\arabic{equation}} 
\setcounter{equation}{0}

Here we prove the following lemma used in section 4 with the notation
introduced there (cf.\ \cite{BFFLS} theorem 5).
\vskip.1cm

\noindent
{\bf Lemma:} {\em If $A(\si,\si')=h(\omega(\si,\si'))$ where $\omega$ is
an antisymmetric 2-form, is an entire analytic function and it satisfies 
(\ref{P2}), (\ref{P3}), then if $A$ is not identically zero, $h(x)$ 
is either $c\sinh\,\mu x$ or $c x$, $\mu,\,c$ being constants.}
\vskip.1cm

\noindent
{\bf Proof:} Putting $\omega(\si,\si')=x$, $\omega(\si',\si'')=y$, 
$\omega(\si'',\si)=z$, (\ref{P2}) becomes
\bez 
h(x-z)h(y)+h(y-x)h(z)+h(z-y)h(x)=0\;.
\eez
Differentiating with respect to $z$ at $z=0$, we find 
\be h'(0)h(y-x)=h(y)h'(x)-h'(y)h(x)\;,\label{Ap1} \ee
where we used that by (\ref{P3}) $h$ is odd, hence $h'$ is even and 
an accent denotes the derivative of $h$ with respect to its argument.
If $h'(0)=0$ then $h'(x)/h(x)=h'(y)/h(y)$ hence $h(x)=c e^{\lambda x}$
and since $h'(0)=0$ either $c=0$ or $\lambda=0$. In both cases $h(x)=0$
since $h$ is odd in $x$. Consequently $h'(0)=\lambda\neq0$. Putting
$\tilde{h}(x):=h(x)/\lambda$, (\ref{Ap1}) becomes
\be 
\tilde{h}(x+y)=\tilde{h}(y)\tilde{h}'(x)+\tilde{h}'(y)\tilde{h}(x)\;.
\label{Ap2}
\ee
Differentiating (\ref{Ap2}) succesively with respect to $x$ and $y$ and
equating the results, we get
\bez \tilde{h}(y)\tilde{h}''(x) = \tilde{h}''(y)\tilde{h}(x)\;. \eez
If $h''(x)$ is not identically zero in any open region, then
\bez  
{\tilde{h}''(x)\over\tilde{h}(x)}={\tilde{h}''(y)\over\tilde{h}(y)}=:\mu
=\mbox{constant}\;,
\eez
everywhere and the result follows (cf.\ footnote 4). 

If $\tilde{h}''(x)=0$ in an open region, then  by (\ref{P3})
$h(x)= c x$ and by the analyticity of $A$ this holds everywhere. \hfill Q.E.D. 
\vskip.1cm

\noindent
{\em Remark:} From (\ref{I7}), we immediately see that $h(x)=x$, or
$h(x)=\sinh\,\mu x$ implies that $[\,,]_A$ is essentially the Poisson
bracket or the Moyal Lie bracket.

\end{document}